\input amstex
\documentstyle{amsppt}
\magnification\magstep1
\input xy
\xyoption{all}
\CompileMatrices

\define\defor{\tilde\twoheadrightarrow}

\define\fib{\twoheadrightarrow}
\define\gat{\overset\sim\to\partial}
\define\dogato{\gat_0*\gat_1}
\define\gatc{\gat_0}
\define\gatu{\gat_1}
\define\gati{\gat_i}
\define\Bh{\hat{B}}
\define\harlot{\Omega A/{\Omega^+A}^\infty}
\define\harlow{\Omega UA/{\Omega^+UA}^\infty}

\define\upa{[U\Cal{PA}]}
\define\pupa{[U\Cal{PA}]^{pol}}
\define\pa{\Cal{PA}}
\define\spa{[\Cal{PA}]}
\define\hpaq{[\Cal{PAQ}]}
\define\paq{\Cal{PAQ}}
\define\poleq{\overset pol\to\equiv}
\define\polho{\overset pol\to\sim}
\define\Ker{\operatorname{Ker}}
\define\holim{\operatorname{holim}}
\define\ili{\underset\leftarrow\to{\holim}}

\define\LF{LF}
\define\eLF{\operatorname\LF}

\define\fiber{fiber}
\define\hofiber{\operatorname\fiber}
\define\rat{\Bbb Q}

\topmatter  %XXXXXXXXXXXX

\title On the derived functor analogy in the Cuntz-Quillen framework
for cyclic homology\endtitle
\author by Guillermo Corti\~ nas\endauthor
\affil Affiliation:  Departamento de
Matem\' atica, Facultad de Ciencias Exactas , Universidad
de La Plata.\endaffil
\leftheadtext{Guillermo Corti\~ nas}
\rightheadtext{Derived Functor Analogy }
\address Guillermo Corti\~ nas,
Departamento
de Matem\' atica, Facultad de Ciencias Exactas, Calles 50 y 115, 
(1900) La Plata, Argentina.\endaddress
\email  willie\@mate.unlp.edu.ar\endemail

\abstract Cuntz and Quillen have shown that for algebras over a field $k$
with $char(k)=0$, periodic cyclic homology may be regarded, in some sense,
as the derived functor of (non-commutative) de Rham (co-)homology. The
purpose of this paper is to formalize this derived functor analogy. We
show that the localization ${De\negthinspace f}^{-1}\pa$ of the category
$\pa$ of countable pro-algebras at the class of (infinitesimal)
deformations exists (in any characteristic) (Theorem 3.2) and that, 
in characteristic zero, periodic cyclic homology is the derived functor
of de Rham cohomology with respect to this localization (Corollary 5.4). We 
also compute the derived functor of rational $K$-theory for algebras over 
$\Bbb Q$, which we show is essentially
the fiber of the Chern character to negative cyclic homology (Theorem 6.2).
\endabstract
\endtopmatter  %XXXXXXXXXXX

\document  
\bigskip
\subhead 0. Introduction\endsubhead
\bigskip
In their paper [CQ2], Cuntz and Quillen show that, if $char(k)=0$,
then periodic cyclic homology may be regarded, in some sense, as the derived
functor of (non-commutative) de Rham (co-)homology. The purpose of this paper
is to formalize this derived functor analogy. We show that the localization
${De\negthinspace f}^{-1}\pa$ of the category $\pa$ of countable pro-algebras 
at the class of (infinitesimal)
deformations exists (in any characteristic) (Theorem 3.2) and that, 
in characteristic zero, periodic cyclic homology is the derived functor
of de Rham cohomology with respect to this localization (Corollary 5.4). We 
also 
compute the derived functor of rational $K$-theory for algebras over $\Bbb Q$, 
which we show is 
essentially
the fiber of the Chern character to negative cyclic homology (Theorem 6.2).
For the construction
of ${De\negthinspace f}^{-1}\pa$, we equip $\pa$ with the analogy of a
closed model category structure, where the analogy of cofibrant objects are 
the quasi-free pro-algebras
and the analogy of trivial fibrations are the deformations. Further, we define
notions of strong and weak nil-homotopy between pro-algebra homomorphisms such
that --as is the case with ``real" model categories ([Q])-- 
${De\negthinspace f}^{-1}\pa$ turns out to be isomorphic to the localization
of $\pa$ at the class of weak nil-homotopy equivalences, and equivalent to
the localization of the subcategory of quasi-free algebras (i.e. the cofibrant
objects) at the class of strong nil-homotopy equivalences (cf. 3.2). 
Of course this 
result would be automatic if the structure we put on $\pa$ 
were a model category (cf [Q]), which we prove it is not (3.6). However 
the analogy we have
is sufficient to prove those localization properties and 
to consider derived functors therefrom. Quillen proves (in [Q]) that a 
functor
between model categories which maps weak equivalences between cofibrant
objects into weak equivalences admits a derived functor. The analogy 
of this result also holds in our setting; it says roughly that if a functor 
$\pa @>>>\Cal C$ remains invariant under pro-power series extensions of 
quasi-free pro-algebras (i.e. $F(A\{ X\} /<X>^\infty)\cong FA$), then its left 
derived functor exists (Theorem 5.2). 
Functors satisfying the latter condition are called Poincar\'e functors,
as the condition that defines them is precisely a Poincar\'e lemma for 
(non commutative) power
series. For example if $F$ satisfies the stronger condition $FA=FA[t]$ then
it is Poincar\'e; such is the case of de Rham cohomology in characteristic
zero. Unless explicitly mentioned, all results in this paper hold over any 
characteristic.
\smallskip
The notion of nil-homotopy used here (although related to) is different
from the usual notion of polynomial (or pol-) homotopy, as used for example in 
Karoubi-Villamayor $K$-theory (see Section 4 below). In fact, a typical
homotopy equivalence under pol-homotopy is the inclusion into the
polynomial pro-algebra 
$B\hookrightarrow B[t]$ which is not an 
equivalence under nil-homotopy. Instead, the inclusion into the power series
pro-algebra $B\hookrightarrow B[t]/<t>^\infty$ is a nil-homotopy equivalence. 
Under nil-homotopy, quasi-free pro-algebras are precisely those having the 
homotopy extension property; other properties of quasi-free pro-algebras 
proven in [CQ1] are shown here to have a natural interpretation in terms of 
homotopy (Theorem 2.1). 
\bigskip
The rest of this paper is organized as follows. In section 1, the notion
of (strong) nil-homotopy is introduced, and its first properties are proved.
Section 2 is devoted to the interpretation of quasi-free pro-algebras
as cofibrant objects with respect to the setting of the previous section 
(Theorem 2.1).
The notion of weak nil-homotopy is introduced in section 3, where the existence
of the localized category ${De\negthinspace f}^{-1}\pa$ is proved 
(Theorem 3.2). 
Section 2
is devoted to the comparison between our notion of nil-homotopy and the 
usual, polynomial homotopy. We prove that the localization at the union
of the classes of nil-deformations and graded deformations exists and can
be calculated as a homotopy category (Theorem 4.1). Section 5 deals with
the formalization of the derived functor analogy of [CQ2]. We establish
sufficient conditions for the existence of left derived functors (Theorem 5.2)
and prove that, in characteristic zero, these conditions are met by the 
de Rham supercomplex functor $A\mapsto XA$ of Cuntz-Quillen (Corollary 5.4).
In section 6 we compute the derived functor of the rational $K$-theory
of rational pro-algebras, (Theorem 6.2) and of the negative cyclic homology 
of pro-algebras over any field (Corollary 6.9).
\bigskip
\definition{Note on Notation} We use most of the notations and notions
established in [CQ 1,2,3]. However, some notations do differ: we write
$\partial_i$ ($i=0,1$) for the natural inclusions $1*0$, $0*1:A @>>> QA=A*A$,
and $qa=\partial_0a-\partial_1 a$. Thus our $qa$ is twice Cuntz-Quillen's.
Also our curvature is minus theirs; here 
$\omega(a,b)=\rho a\rho b-\rho (ab)$. In this paper, the superscript
$B^+$ on a graded algebra $B$ denotes the terms of positive degree, and
not the even degree part as in {\it op. cit.}. The even and odd terms
are indicated by $B^{even}$ and $B^{odd}$. If $A$ is a pro-algebra indexed
by $\Bbb N$, then the map $A_{n+1} @>>> A_n$ is referred to as the structure
map and is named $\sigma$ or $\tau$ (subscripts are mostly omitted). Since 
for the most part we make no assumptions on $chark$, none of the results of 
{\it op. cit.} which involve dividing by arbitrary integers holds. Such is the 
case of the
isomorphism between $QA$ and the de Rham algebra with Fedosov product ([CQ1]),
-- as it assumes $2\ne 0$-- which we do not use. We do use the fact that 
$qA^n/qA^{n+1}\cong\Omega^nA$ as $A$-bimodules, which does hold even if $2$ 
is not invertible. On the other hand the isomorphism between the tensor
algebra $TA$ and the algebra of even differential forms holds in any
characteristic with the same proof as in [CQ1].\enddefinition
\bigskip
\bigskip
\bigskip
\subhead 1. A Closed Model Category Analogy\endsubhead
\smallskip
\item{1.0}\ \ We consider associative, non-necessarily unital algebras
over a fixed ground field $k$. We write $\Cal A$ and $\Cal V$ for the
categories of algebras and vector spaces and $\pa$ and $\Cal{PV}$ for
the corresponding pro-categories. As in [CQ3] we consider only countably
indexed pro-objects. A map $f\in \Cal{PA}(A,B)$ is called a {\it fibration}
if it admits a right inverse as a map of pro-vector spaces, i.e. there
exists $s\in {\Cal{PV}}(B,A)$ such that $fs=1$. Fibrations are denoted
by a double headed arrow $\twoheadrightarrow$. By a (nil-) deformation
($\defor$) of a pro-algebra $A$ we mean a fibration onto $A$ which is 
isomorphic to one
of the form $P/K^\infty\defor P/K$. Equivalently, $p:B\defor A$ is
a deformation iff it is a fibration and for $K=\Ker (p)$ we have $K^\infty =0$.
For example the map:
$$
UA:=TA/JA^\infty\overset\pi^A\to\defor A 
$$
is a deformation, and is initial among all deformations with values in $A$.
That is if $p:B\defor A$ is a deformation then there exists a map $f:UA\longrightarrow B$
with $pf=\pi^A$. In particular if $A$ is quasi-free in the sense of [CQ3]
then $p$ is split in $\pa$ (because $\pi^A$ is). Deformations 
admitting a right inverse shall be called {\it deformation retractions}; thus
$A$ is quasi-free iff every deformation $B\defor A$ is a retraction (or
$A$ is a retract of every deformation onto it). It follows that
quasi-free pro-algebras are precisely those pro-algebras $A$ such that the
map $0\rightarrowtail A$ has the left lifting property (LLP) with respect to
deformations. Thus we have the analogy of closed model category ([Q]) where
fibrations are as above, trivial fibrations are deformations, and cofibrant
objects are quasi-free algebras. To pursue this analogy a step further,
we define our weak nil-equivalences (or wne's) as follows. We say that a map $f\in\pa$
is a wne if any functor defined on $\pa$ and taking values in some 
category $\Cal C$ which inverts (i.e. maps to isomorphisms) all 
nil-deformations also inverts $f$. Functors which invert wne's are called
nil-invariant. We shall show that the localization of $\pa$ with
respect to deformations exists, whence $f$ is a wne iff it is inverted
upon localizing. For completeness, we call a map $f$ quasi-free if it
has the LLP with respect to deformations. Thus quasi-free maps play the
r\^ole of cofibrations. I hurry to point out that the above notions of fibration,
cofibration, and weak equivalence DO NOT make $\pa$ into a closed model
or even into a model category. Indeed, if the map $0\longrightarrow A$ factors
as a weak equivalence followed by a fibration then $A$ is weak equivalent
to $0$ (3.5). As there are pro-algebras which are not equivalent to zero, 
axiom M2 for a model category ([Q]) does not hold. The latter problem would 
be solved if we allowed free maps of the form $A\longrightarrow A*TV$ to be weak equivalences;
in fact any map $A\longrightarrow B$ factors as $A\longrightarrow A*TB$ followed
by $a\mapsto f(a), \rho b\mapsto b$. This simply means that there are nil-invariant
functors which do not invert free maps. 
\smallskip
The notion of weak equivalence defined above may be expressed as the weak
homotopy relation associated to a notion of strong homotopy 
between pro-algebra homomorphisms. The definition of this strong homotopy
is the subject of the next subsection.
\bigskip
\subhead{\smc Cylinders and nil-homotopy 1.1}\endsubhead
\medskip
The {\it cylinder} of a pro-algebra $A$ is the following pro-algebra:
$$
Cyl(A):=QA/qA^\infty\tag1
$$
Here $QA=A*A$ is the free product (or coproduct, or sum) and 
$qA=\Ker (QA\longrightarrow A)$ is
the kernel of the folding map. We write $\partial_0=1*0$ and 
$\partial_1=0*1$ for the canonical inclusions $A\longrightarrow QA$,\ \
$\overset\sim\to\partial_0*\overset\sim\to\partial_1:QA\longrightarrow CylA$
for the completion map, and $p=p_A:CylA\defor A$ for the the completion
of the folding map $\mu:QA\longrightarrow A$. We have a commutative diagram:
$$
\xymatrix{
QA\ar[d]_\mu\ar[dr]^\dogato&\\
A &\ar[l]_\sim CylA}\tag2
$$
One checks that $\overset\sim\to\partial_0*\overset\sim\to\partial_1$ is
quasi-free if A is, whence $CylA$ is a cylinder object in the sense
of [Q, 1.5. Def. 4]. Given homomorphisms $f,g:A\longrightarrow B$, we
write $f\equiv g$ if there exists a map $h:CylA\longrightarrow B$ making
the following diagram commute:
$$
\xymatrix{QA\ar[r]^{f*g}\ar[d]_\dogato &B\\
CylA\ar[ur]_h&}\tag3
$$
Note that as $QA\longrightarrow CylA$ is an epimorphism (although not
a fibration), if a homotopy (i.e. a factorization through $CylA$ ) exists, 
it must be unique. For example if $A$ and $B$ are algebras, then $f\equiv g$ iff there
exists n such that for all $a_1,\dots , a_n\in A$, we have
$$
(f(a_1)-g(a_1))\dots (f(a_n)-g(a_n))=0
$$
and the homotopy is the map sending the class of $qa$ to $f(a)-g(a)$.
One checks that $\equiv$ is a reflexive and symmetric relation, and that
it is compatible with composition on the left:
 $f_0\equiv f_1$ $\Rightarrow$ $f_2f_0\equiv f_2f_1$ (whenever composition
makes sense).
It follows that the equivalence relation $\sim$ generated by $\equiv$ is compatible with composition on both
sides. We say that $f$ and $g$ are (nil-) homotopic if $f\sim g$. We write
$[\pa]$ for the category having the same objects as $\pa$ and
as morphisms the sets of equivalence classes:
$$
[A,B]:=\Cal{PA}(A,B)/\sim
$$
A map $f\in \pa$ is called a {\it strong} nil-homotopy equivalence
if its class is an isomorphism in $[\Cal{PA}]$.
\bigskip
\remark{Remark 1.2} The homotopy relation defined above may also be
defined in terms of $n$-{\it fold cylinders}. Set $Cyl^1A:=CylA$, 
$\gati^1=\gati$ and define the $n$-fold cylinder inductively by
the pushout diagram:
$$
\CD
  A         @>\gatc>>      Cyl^1A\\
@V\gatu^{n-1}VV              @VVV\\
Cyl^{n-1}A  @>>>              Cyl^nA\\ 
\endCD
$$
Define $\gatc^n$ as the composite map $A@>\gatc^{n-1}>>Cyl^{n-1}A@>>>Cyl^nA$
and $\gatu^n$ as the composite $A@>\gatu^1>>CylA @>>>Cyl^nA$. One checks that two maps 
$f,g:A\longrightarrow B$ are homotopic iff there exist $n$ and $h:Cyl^nA\longrightarrow B$
such that the following diagram commutes:
$$
\xymatrix{QA\ar[d]_{\gatc^n*\gatu^n}\ar[r]^{f*g}& B\\
Cyl^nA\ar[ur]_h}
$$
The map $h$ in the diagram above will be called a homotopy between $f$ and 
$g$.\endremark

\bigskip
The following lemma establishes a relation between the nil-homotopy equivalences
just defined and the weak nil-equivalences of 1.0. above.
\bigskip
\proclaim{Lemma 1.3} Let $f:A\defor B$ be a deformation retraction.
Then $f$ is a strong nil-homotopy equivalence.\endproclaim
\demo{Proof} We have to prove that $g=sf\sim 1$. Upon re-indexing,
we can assume $f=\{ f_n:A_n\longrightarrow B_n\}$, $s=\{s_n:B_n\longrightarrow A_n\}$
are inverse systems of maps commuting with the structure maps
$\sigma=\sigma_n$, that $\sigma f_ns_n=\sigma$ and that for $K_n=\Ker f_n$
we have $K_n^n=0$. Then for $a\in A_n$, we have $f(\sigma (g_n*1)qa)=\sigma (fsfa-fa)=0$,
from which $\sigma (g_n*1(qa))\in K_{n-1}$. Thus $\sigma (g_n*1)(qA_n)^n=0$ whence 
$g*1:QA\longrightarrow B$ factors through $CylA$, and $g\equiv 1$.\qed\enddemo
\bigskip
\subhead{2. Quasi-free Algebras and the Homotopy Extension Property}\endsubhead
\bigskip
An interesting feature of nil-homotopy is that quasi-free algebras are precisely
those having the homotopy extension property with respect to deformations.
This fact is proven in Theorem 2.1 below. First we need:
\bigskip
\subhead{\smc Power pro-algebras, power spans and power deformations 2.0}\endsubhead
By a {\it graded} pro-algebra we mean a non-negatively graded object in
$\pa$, i.e. a pro-algebra $B$ together with
a direct sum decomposition of pro-vector spaces: $B=\bigoplus_{n=0}^\infty B^n$
such that the multiplication map $B\otimes B=\bigoplus B^n\otimes B^m\longrightarrow B$
maps $B^n\otimes B^m$ into $B^{n+m}$. Thus $B^+=\bigoplus_{n=1}^\infty B^n$
is a two-sided ideal in $B$, in the sense that multiplication maps $B^+\otimes B$
and $B\otimes B^+$ into $B^+$. It is straightforward to show that every
graded pro-algebra is isomorphic--by a homogeneous isomorphism-- to an inverse
system of graded algebras and homogeneous maps. The {\it power} pro-algebra associated with
$B$ is the pro-algebra $\hat B:=B/B^{+^\infty}$. Thus a power pro-algebra 
is a particular kind of graded algebra. For instance if $A$ is an algebra 
then the power pro-algebra associated to the polynomials in a set $X$ is
the pro-algebra $\{ A\{X\}/<X>^n\}$, whose completion is the power series
algebra in the non-commutative variables $X$. More generally, one considers
the tensor algebra 
$T_{\tilde A}(M)=T_0(A)\bigoplus T^1(A)\bigoplus T^2(A)\dots =
A\oplus M\oplus M\otimes_{\tilde A} M\oplus\dots$
whose associated power algebra is $\hat T_{\tilde A}(M)=
\{\bigoplus_{i=1}^n T^i(M): n\in\Bbb N\}$
and when $M$ is the free module on a set $X$ one recovers the polynomial and
power series algebras. These constructions can be copied for pro-algebras, 
pro-sets and
pro-modules with the obvious definitions. However in general the free pro-module
associated with a pro-set is not proyective, as it doesn't have the LLP with
respect to all epimorphisms, but only with respect to fibrations (cf.[CQ3]).
We use the following special notations. If $V$ is a pro-vector space and
$I\vartriangleleft A$ is an ideal in a pro-algebra, we write $P_A(V)$ for the power algebra associated
with $T_{\tilde A}(\tilde A\otimes V\otimes\tilde A)$ and $G_I(A)$ and $\hat G_I(A)$
for the graded pro-algebra $A\oplus I/I^2\oplus I^2/I^3\oplus\dots$ and
its associated power algebra. If $B$ is a graded pro-algebra and $u:A\longrightarrow B^0$
is a homomorphism, then by a {\it power span} of $u$ we mean a $k$-linear map 
$T=\sum_{n=1}^\infty D_n: A\longrightarrow \Bh^+$
such that the following diagram commutes:
$$
\CD
A\otimes A                             @>\text{multiplication}>>          A\\
@V u\otimes T+T\otimes u-T\otimes T VV                                @VVTV\\
\Bh\otimes \Bh\oplus \Bh\otimes \Bh\oplus  \Bh\otimes \Bh @>\text{sum+multiplication}>> \Bh^+\\
\endCD\tag4
$$
Briefly, we write 
$$
T(xy)=uxTy+Txuy-TxTy\tag4'
$$
to indicate the diagram above --even
if $A$ and $B$ are not algebras. For example the ordinary Taylor span:
$$
k[x]\longrightarrow k[x][[y]]=\{ k[x,y]/<y>^n\}, f(x)\mapsto \{\sum_{i=0}^n\frac{f^{(i)}(y)}{i!}\}
$$
is a power span of the canonical inclusion. Note that the image of $f(x)$ in
$k[x,y]/<y>^n$ is just the class of $f(x)-f(y)$ and is therefore defined
in any characteristic; if $f(x)=\sum_{i=0}^n a_ix^i$ then $\frac{f^{(i)}(y)}{i!}$ is just
short for $\sum_{j=0}^{n-i} \binom{n}j a_{i+j}$ which is defined everywhere.
Note also that any power span $T$ induces a homomorphism
$h:CylA\longrightarrow \Bh$ with $h\gatc =u$, which is a homotopy between $u$ and
$h\gatu$. Conversely if $h$ is a homotopy starting at $u$, then 
$T:A\overset q\to\longrightarrow qA\longrightarrow qA/qA^\infty\overset h \to \longrightarrow \Bh$
is a power span. Thus a power span is a special kind of homotopy where the target
is a power algebra. By an $n$-{\it truncated span} we mean a linear map $T_n:A\longrightarrow B/{B^+}^{n+1}$
satisfying (4'). For example if $T$ is a power span then $T_n:A\overset T\to\longrightarrow B/{B^+}^\infty\longrightarrow B/{B^+}^{n+1}$
is an $n$-truncated power span. Finally, by a {\it power deformation retraction} 
we mean a deformation retraction of the form $\Bh\rightarrow B_0$ where $B$ is
a graded algebra.
\bigskip
\proclaim{Theorem 2.1}{\rm (Compare [CQ1]).} The following conditions
are equivalent for a pro-algebra $A$.
\medskip
\roster
\item"{(i)}"({\rm LLP}) $A$ is quasi-free.
\medskip

\item"{(ii)}"({\rm Power Span Extension}) If $B$ is a graded algebra and $u:A\longrightarrow B^0$ is a homomorphism
then any truncated span $T_n:A\longrightarrow B/{B^+}^{n+1}$ lifts to a power
span $T:A\longrightarrow \Bh$.
\medskip
\item"{(iii)}"({\rm Tubular Neighborhood}) If $f:B\defor A$ is a 
deformation with kernel $I$ and $B$ is 
quasi-free, then there is an isomorphism 
$\iota:B\overset\cong\to\rightarrow\hat{G}_I(B)$ such that $f\iota$ is
the projection $\hat{G}_I(B)\defor B/I=A$.
\medskip
\item"{(iv)}"({\rm Even Forms}) There is a pro-algebra isomorphism
 $UA\cong \Omega^{even} A/$\linebreak ${\Omega^{even+}}^{\infty} A$
which makes the following diagram commute:
$$
\xymatrix{
UA\ar[d]_{\pi^A\wr}\ar[r]^\simeq &\Omega^{even} A/{\Omega^{even+}}^\infty
A\ar[dl]\\
A&}
$$
\medskip
\item"{(v)}"({\rm de Rham Algebra}) There is a pro-algebra isomorphism 
$CylA\cong \Omega A/{\Omega^+}^\infty A$
which makes the following diagram commute:
$$
\CD
CylA     @>\simeq>>   \Omega A/{\Omega^+}^\infty A\\
@VVV                        @VVV\\
A\oplus qA/qA^2   @>\simeq>> A\oplus\Omega^1A\\
\endCD
$$
Here the bottom arrow is the canonical isomorphism $aqb\mapsto adb$.
\medskip
\item"{(vi)}" {\rm (Homotopy Extension)} Given any commutative solid
arrow diagram:
$$
\xymatrix{A\ar[r]\ar[d]_{\gatc} &  B\ar[d]^{\wr f}\\
CylA\ar@{..>}[ur]\ar[r]& C}
$$
where $f$ is a deformation, the dotted arrow exists and makes it commute.
\endroster\endproclaim
\demo{Proof} (i)$\Rightarrow$(ii): Write $T_n=\sum_{i=1}^n D_i$ where $D_i$                                                             
is the part of degree $i$; also let $D_0=u$. Thus $u_n=u+T_n=\sum_{i=0}^n D_i$
is a homomorphism, from which the following identity follows:
$$
-\delta D_i=\sum_{j=1}^i D_j\cup D_{i-j}\qquad\qquad (0\le i\le n)\tag5
$$
Here the maps $D_i$ are regarded as $1$-cochains with values in $B$, the
cup product is the composite of $D_j\otimes D_{i-j}$ with the multiplication
map $B\otimes B\longrightarrow B$ and $\delta$ is the Hochschild co-boundary
map---as defined by the appropriate arrow diagram. We must prove
that a $k$-linear map $D_{n+1}:A\longrightarrow B^{n+1}$ exists so that
$$ 
-\delta D_{n+1}=\sum_{i=1}^n D_i\cup D_{n+1-i}\tag5'
$$
holds. It is straightforward to check that the right hand side of (5') is actually
a cocycle, whence also a coboundary, as $A$ is quasi-free. Explicitly, if
$g:\Omega^2(A)\longrightarrow B_{n+1}$ is the bimodule homomorphism induced by the
right hand side of (5') and if
$f:A\longrightarrow \Omega^2(A)$ satisfies $-\delta f=d\cup d$, then we can
take $D_{n+1}=gf$.
\smallskip
(i)$\iff$(iii): If (iii) holds then $UA\defor A$ is a retraction, whence
$A$ is quasi-free. Supose conversely that (i) holds. Because $A$ is 
quasi-free, we have direct sum
decompositions $B=A\oplus I$, and 
$B/I^2=A\oplus I/I^2=\hat{G}_I(B)/{\hat{G}_I(B)^+}^2$. Write 
$u:B\defor A\hookrightarrow \hat{G}_I(B)$ for the composite map,
and $p_1:B\defor \hat{G}_I(B)/{\hat{G}_I(B)^+}^2$ for the projection. Because $B$ is
quasifree, the truncated span 
$T_1=p_1-u:B\rightarrow {\hat{G}_I(B)^+}/{\hat{G}_I(B)^+}^2$ extends to
a power span $T:A\rightarrow {\hat{G}_I(B)^+}/{\hat{G}_I(B)^+}^{\infty}$
(by (ii)). It is clear that $T$ induces the identity on $I/I^2$;
further, one checks --using (5)-- that it also induces the identity on
$I^n/I^{n+1}$. It follows that $p:u+T$ is an isomorphism.
\smallskip
(iii)$\Rightarrow$(iv): Applying (iii) to $\pi^A:UA\defor A$, we get 
\smallskip
 $UA\cong\hat{G}_{JA/JA^{\infty}}(UA)=\Omega^{even} A/{\Omega^{even+}}^{\infty} A$.
\smallskip
(iv)$\Rightarrow$(i): Analogous to (iii)$\Rightarrow$(i).
\smallskip
(ii)$\Rightarrow$(v): By (ii), we can lift the de Rham derivation
$d:A\longrightarrow \Omega^1A$ to a power span $T:A\longrightarrow \Omega A/{\Omega^+}^\infty A$
of the identity map $A=\Omega^0A$. By the discussion above, $1+T$ induces
a homomorphism $h:CylA\longrightarrow\Omega A/{\Omega^+}^\infty A$ such that $hq=T$.
In particular, $h$ induces the canonical $A$-bimodule isomorphism 
$qA/qA^2\cong\Omega^1A$ mapping $q$ to $d$. Thus we have $hq=d+D$, where
$D(A)\subset \Omega^{\ge 2}/{\Omega^{\ge 2}}^\infty$. It follows that the composite
$A^{\otimes n}\overset hq^{\otimes n}\to\longrightarrow\Omega^+/{\Omega^+}^\infty\longrightarrow\Omega^+/{\Omega^+}^{n+1}$
is just the cocycle $d^{\cup n}$, whence the induced bimodule homomorphism
$qA^n/qA^{n+1}\cong\Omega^nA$ is the canonical isomorphism, and the proof
ensues.
\smallskip
(v)$\Rightarrow$(i): By virtue of (5), if $T_2=d+D_2:A\longrightarrow\Omega^1A\oplus\Omega^2A$
is the $2$-span induced by $\gatu$, then $-\delta D_2=d\cup d$, whence $A$ is 
quasi-free.
\smallskip
(vi)$\Rightarrow$(v): Since $\harlot\defor\Omega^0A\oplus\Omega^1A$
is a deformation, there exists a homomotopy $h:CylA\longrightarrow \Omega A/{\Omega^+}^\infty A$
lifting the homotopy $1\equiv 1+d$. The same argument as in the proof of 
(ii)$\Rightarrow$(v) shows that $h$ is an isomorphism.
\smallskip
(i)$\Rightarrow$(vi): As $0\longrightarrow A$ is quasi-free, so are $\partial_0$
and $\gatc$.\qed\enddemo
\bigskip
\example{Example 2.2} Let $A$ be an algebra, and let $UA=TA/JA^\infty$ its 
universal quasi-free model. By the theorem above, we have $CylUA\cong\harlow$.
We want to give an explicit isomorphism $CylUA\cong\harlow$
as well as to show that in this particular case, we also have an isomorphism
$$
\harlow\cong P_{UA}(A)
$$
First of all, we observe
that given a vector space $V$, we have isomorphisms:
$$
\split
QTV&\cong T(V\oplus V)\cong T(V\oplus qV)\cong T(V)*T(qV)\\
   &\cong T_{\widetilde{TV}}(\widetilde{TV}\otimes V\otimes\widetilde{TV})\cong \Omega TV\\
\endsplit
$$
Here $qV=\{(v,-v): v\in V\}$ and the isomorphism $V\oplus V\cong V\oplus qV$
maps $(v,0)=\partial_0v$ to itself while $\partial_1v\mapsto qv$. Thus the composite
isomorphism $\alpha:QTV\vec{\cong}\Omega TV$ maps $qv$ to $dv$ and $\partial_0x$ to $x$
($v\in V, x\in TV$). In particular this holds when $V=A$; in this case
$\alpha$ maps the
ideal $<JA>\subset QTA$ generated by $JA$ (which we identify with its image
through $\gatc$) into the ideal $<JA>\subset \Omega TA$,
and $qTA$ into $\Omega^+ TA$. It follows that $\alpha$ induces an isomorphism
$QTA/{\Cal F}^\infty\cong\Omega TA/{\Cal G}^\infty$, where $\Cal F$ and
$\Cal G$ are respectively the $<JA>+qTA$ and $<JA>+\Omega^+TA$-adic 
filtrations. On the other hand we have $CylUA=QTA/{\Cal F'}^\infty$ and 
$\harlow =\Omega TA/{\Cal G'}$ where ${\Cal F'}=<JA^n>+<q(JA^n)>+(qTA)^n$
and ${\Cal G'}=<JA^n>+<dJA^n>+(\Omega^+TA)^n$. 
We have inclusions:
\smallskip 
${\Cal F}^n\supset {\Cal F'}^n\supset {\Cal F"}^n=<JA^n>+(qTA)^n$
\smallskip
and
\smallskip
${\Cal G}^n\supset {\Cal G'}^n\supset {\Cal G"}^n=<JA^n>+(\Omega^+TA)^n$
\smallskip
Lemma 2.3. below shows that for $N$ sufficiently large, we also have inclusions
${\Cal F"}^n\supset {\Cal F}^N$ and ${\Cal G"}^n\supset {\Cal G}^N$. It follows
that $\alpha$ induces the isomorphism $CylUA\vec{\cong}\harlow$ and that 
$\harlow=\Omega TA/G"^\infty=P_{UA}(A)$\endexample
\bigskip
\proclaim{Lemma 2.3} Let $A\subset B$ be algebras and let $\epsilon:B\rightarrow A$
be a homomorphism such that $\epsilon a=a, (a\in A)$. Set $I=\Ker\epsilon$,
and let $J\subset A$ be an ideal. Consider the following filtration in
$B$:
$$
B\supset {\Cal F}^n=<J^n>+I^n
$$
Then there is an isomorphism:
$$
B/{\Cal F}^\infty\cong B/(<J>+I)^\infty
$$
\endproclaim
\smallskip
\demo{Proof} Let ${\Cal G}^n=<J>^n+I^n$. It is straightforward to check
that $(<J>+I)^{2n}\subset{\Cal G}^n$, whence $B/(<J>+I)^\infty\cong B/{\Cal G}^\infty$.
Thus we must prove that $B/{\Cal G}^\infty\cong B/{\Cal F}^\infty$. It is clear
that ${\Cal G}^n\supset {\Cal F}^n$. I claim that for $N=n^2+n-1$, we also have
${\Cal G}^N\supset {\Cal F}^n$. To prove the claim --and the lemma-- it
suffices to show that $<J>^N\subset{\Cal F}^n$. Every element of $<J>^N$ is
a sum of products of the form:
$$
(j_1+i_1)\dots(j_N+i_N)\qquad\qquad (j_r\in J, i_r\in I)
$$
After fully expanding the product above, we get a large sum in which those
terms not in $I^n$ have at most n-1 $i$'s and at least $n^2$ $j$'s. Therefore,
in each such term, 
at least $n$ of the $j$'s must appear side by side, forming a string. Hence
the term in question lives in $<J^n>$.\qed\enddemo
\bigskip
\remark{Remark 2.4} The de Rham pro-algebra $\harlot=\{\bigoplus_{r=0}^n\Omega^rA_{n+1}\}$, 
of a pro-algebra $A=\{ A_n\}$, together with the
natural differentials $b$ and $d$ and the Karoubi operator $\kappa$, can
be regarded as a pro-{\it truncated} mixed $DGA$ in the sense of [Kar].
Indeed, the identity:
$$
bd\omega +db\omega =\omega -\kappa\omega
$$
holds in $\Omega^r(A_{n+1})$ for $r<n$ and in 
$\Omega_{\natural}^nA_{n+1}=\Omega^nA_{n+1}/[\Omega^0A,\Omega^nA]$ for $r=n$.
Thus:
$$
\theta\Omega(A)=((\oplus_{r=0}^{n-1}\Omega^rA_{n+1})\oplus\Omega_{\natural}^nA_{n+1}, B+b)
$$
is a pro-differential graded vector space, equipped with an even-odd gradation.
This is the pro-complex of [CQ-2]; if $k\supset \rat$, it is homotopy
equivalent to the (short) de Rham pro-complex:
$$
XUA:\Omega^0UA\overset{\overset{\underset{\natural}\to{\, d}}\to\longrightarrow}\to{\underset b\to\longleftarrow}\Omega_{\natural}^1UA
$$
In any characteristic, we still have $\theta\Omega UA\approx XUA$ for every
algebra $A$ and $\theta\Omega R\approx XR$ for every quasi-free algebra $R$.
In particular $CylR$ carries all the relevant information
for the cyclic homology of $R$.
\endremark
\bigskip
\subhead{3. The Homotopy Category}\endsubhead 
\bigskip
\subhead{\smc Weak nil-homotopy 3.0}\endsubhead
We write $[U\Cal{PA}]$ for the category having the same objects as 
$\Cal{PA}$ and where the set of maps from $A$ to $B$ is $[UA,UB]$. We 
have a functor $\gamma:\Cal{PA}\rightarrow\upa$, $A\mapsto A$, $f\mapsto [Uf]$.
Two maps $f,g\in \Cal{PA}(A,B)$ shall be called {\it weakly nil homotopic}
if $\gamma f=\gamma g$; by a weak nil homotopy equivalence we shall mean
a map $f\in \Cal{PA}$ such that $\gamma f$ is an isomorphism. We show
below that the class of weak nil homotopy equivalences is precisely
the class of weak nil equivalences as defined in 1.0 above, and that
$\gamma$ is the localization of $\pa$ at this class. Further,
we show that $\upa$ is equivalent to the strong homotopy category $\hpaq$ of
quasi-free algebras. First we need:
\bigskip
\proclaim{Lemma 3.1} The functor $U:\Cal{PA}\rightarrow {\paq}$ 
carries fibrations to fibrations and deformations to deformations.\endproclaim
\demo{Proof} Let $f=\{f_n:A_n\fib B_n\}$ be a fibration, and let $t={t_n:B_n\rightarrow A_n}$
be a section of $f$ in $\Cal{PV}$. Upon re-indexing, we can assume that 
$ft\tau=\tau$ for the structure map of $B$. We want to construct a linear
section $\hat t$ of $Uf$ lifting $t$. Consider the following
composite of linear maps:
$$
s_n:\frac{TB_n}{JB_n^n}\overset\sim\to\rightarrow\bigoplus_{i=0}^{n-1}\Omega^{2i}B_n\hookrightarrow\bigoplus_{i=0}^{\infty}\Omega^{2i}B_n\cong TB_n
$$
Note that $s_n$ is a linear section of $TB_n\fib TB_n/JB^n$. Consider
the composite $\hat t_n:TB_n/JB^n @>s_n>>TB_n@>Tt_n>>TA_n@>>>TA_n/JA_n^n$;
then $\hat t_n$ commutes with $\tau$ and 
$$
\multline
\hat t_n(\rho b_0\omega (b_1,b_2)\dots\omega(b_{2l-1},b_{2l})=\\
=\rho t_n(b_0)(\omega (t_nb_1,t_nb_2)+\rho\omega_{t_n} (b_1,b_2))\dots 
(\omega(t_nb_{2l-1},t_nb_{2l})+\rho\omega_{t_n} (b_{2l-1},b_{2l}))\\
\endmultline
$$
for $0\le l\le n-1$. Here $\rho:A\rightarrow TA$ is the canonical section,
$\omega(a,b)=a\otimes b-ab$ is the curvature of $\rho$ and $\omega_{t_n}$ is 
the curvature
of $t_n$. Now since $ft\tau=\tau$, we have $\omega_{t_n}(b,b')\in \Ker\tau^B 
f_n$  ($b, b'\in B_n$)
and $\rho\omega_{t_n}(b,b')\in \Ker\tau^{TB}Tf_n$. It follows that $Uf_n\hat t_n\tau_n^{UB}=\tau_n^{UB}$,
whence $Uf$ is a fibration. Suppose further that $f$ is also a deformation,
and let $K=\Ker f$; we can assume $K_n^n=0$. Let $L=\Ker Uf$; if $l\in L_n^n$ then
$\pi_n^A l\in K_n^n=0$, hence $L_n^n\subset JA_n/JA_n^n$, and $L_n^{n^2}=0$.\qed\enddemo
\bigskip
\proclaim{Theorem 3.2} (Compare [Qui, 1.13, Th.1])
\roster
\item"{(i)}" Strong nil-homotopy equivalences are precisely those
maps which are inverted by every functor which inverts {\rm deformation retractions}.
Weak nil-homotopy equivalences are precisely those maps in $\pa$
that are inverted by every nil-invariant functor, i.e. every functor which 
inverts {\rm all deformations}.
\medskip
\item"{(ii)}" The functor $\pa\rightarrow [\pa]$ is the localization of
$\pa$ at the class of {\rm deformation retractions}, the functor 
$\paq\rightarrow\hpaq$
is the localization at the class of {\rm power deformation retractions}, and
the functor $\gamma:\pa\rightarrow\upa$ is the localization at the class of 
{\rm all deformations}. There is a category equivalence: $\upa\approx\hpaq$.
\endroster
\endproclaim
\demo{Proof} (i) Let $se$ be the class of maps inverted by every functor
which inverts deformation retractions and let $se'$ be the class of strong
homotopy equivalences. By virtue of Lemma 1.3, the functor $\pa\rightarrow [\pa]$
inverts deformation retractions, whence $se\subset se'$. Conversely, 
if $F$ inverts deformation retractions then it inverts
$CylA\defor A$, and also $\gati$, $i=0,1$. Thus $F$ maps congruent
maps to the same map; further, since $f\overset F\to\sim g\iff Ff=Fg$ is
an equivalence relation, $F$ also maps nil-homotopic maps to the same map,
and strong nil-equivalences to isomorphisms. This proves the first assertion
of (i). Next, write $\omega$ and $\omega'$ for the classes 
of weak nil-equivalences (as defined in 1.0 above) and weak nil-homotopy 
equivalences. We have to prove that $\omega=\omega'$. In view of
Lemmas 1.3 and 3.1, the functor $\gamma$ is nil-invariant, whence 
$\omega\subset\omega'$. Now let $F:\pa\rightarrow\Cal C$ be a nil-invariant
functor, and let $f\in\omega'(A,B)$. Because $F\pi^A$ and $F\pi^B$ are 
isomorphisms in $\Cal C$, $Ff$ will be an isomorphism iff $FUf$ is. By definition,
the fact that $f\in\omega'$ means that $Uf$ is a strong equivalence, and
therefore is inverted by $F$. Thus $\omega=\omega'$.
\smallskip
(ii) The first assertion of (ii) is immediate from the proof of the first assertion of (i). 
The second assertion follows similarly, in view of 2.1-iii).
Now let $F$ be a nil invariant functor as above.
We have to show that $F$ factors as $F=\tilde{F}\gamma$ for some 
$\tilde F:\upa\rightarrow\Cal C$, and that such $\tilde F$ is unique.
We put $\tilde F(A)=F(A)$ and for $[f]\in\upa(A,B)$, we set 
$\tilde F[f]=F\pi^B Ff (F\pi^A)^{-1}$. It is clear that $\tilde F$ is 
well-defined and that $F=\tilde F\gamma$. Now suppose $G$ is another functor
with the same property as $\tilde F$. Then $GA=A$ on objects and if $f\in\pa (A,B)$
then $G$ must map $[Uf]$ onto $FUf=\tilde F [Uf]$. Since any map $[g]\in\upa(UA,UB)$
factors as $[g]=[\pi^{UB}][Ug][\pi^{UA}]^{-1}$, it suffices to prove
that $[\pi^{UA}]=[U\pi^A]$. But both $\pi^{UA}$ and $U\pi^A$ are left inverse
to the same map $\iota:UA\rightarrow U^2A$ induced by $T\rho:TA\rightarrow T^2A$,
whence (by Lemma 1.3) $[\pi^{UA}]=[\iota]^{-1}=[U\pi^A]$. This proves the
third assertion. By the proof of (i), the functor 
$\gamma:\paq\rightarrow\upa$ induces
a functor $\overline{\gamma}:\hpaq\rightarrow\upa$. Let $\gamma':\upa\rightarrow\hpaq$,
$A\mapsto UA$, $[f]\mapsto [f]$. Then $[\pi^R]:\gamma'\overline{\gamma}(R)=UR\rightarrow R$
and $[\pi^{UA}]:\overline{\gamma}\gamma'(A)=UA\rightarrow A$ are natural
isomorphisms $\overline{\gamma}\gamma'\overset\cong\to\rightarrow 1$ and
$\gamma'\overline{\gamma}\overset\cong\to\rightarrow 1$. This concludes the
proof.\qed\enddemo
\bigskip
\proclaim{Corollary 3.3} Let $f,g:A\rightarrow B$ be pro-algebra
homomorphisms. We have:
\roster
\item"{(i)}" Strong $\Rightarrow$ Weak: If $f$ is a strong equivalence
then it is also a weak equivalence. If $f$ and $g$ are strongly
nil-homotopic then they are also weakly homotopic.
\medskip
\item"{(ii)}" Weak $\Rightarrow$ Strong: The converse of (i) holds
if $A$ and $B$ are quasi-free.
\endroster
\endproclaim
\demo{Proof} As $CylA\rightarrow A$ is a deformation, any nil invariant
functor maps strong equivalences into isomorphisms and homotopic maps
to the same map. In particular, this happens with the localization
functor $\gamma$, proving (i). Part (ii) follows from the identities:
$[A,B]=\hpaq (A,B)=\upa (A,B)=[UA,UB]$.\qed\enddemo
\medskip
By defintion, the class $De\negthinspace f$ of deformations sits into the intersection
of the class $we$ of weak equivalences and the class $Fib$ of fibrations.
The proposition below shows that in fact $De\negthinspace f=we\cap Fib$. In particular
this proves that quasi-free maps are precisely those having the LLP with
respect to those fibrations which are weak equivalences.
\bigskip
\proclaim{Proposition 3.4} A fibration is a deformation iff it is
a weak equivalence.\endproclaim
\demo{Proof} If $f$ is deformation then it is a weak equivalence 
by definition of the latter. Suppose now
$f:A\fib B$ is a fibration
and a weak equivalence, and write $K=\Ker f$.
Upon re-indexing, we can assume $f$ is an inverse system of epimorphisms
$\{f_n:A_n\fib B_n\}$ commuting with structure maps.
We must prove $K^\infty=0$. I claim it suffices to check this for the particular
case when $f$ is a strong equivalence. For if $f$ is a weak equivalence and
a fibration then $Uf$ is both a strong equivalence (by 3.3) and a fibration
(by 3.1). Whence, if we know the proposition for strong equivalences, 
we have $\Ker Uf^\infty=0$. Now a little diagram chasing shows that $\Ker Uf_n\fib K_n$
is an epimorphism ($n\ge 1$), whence also $K^\infty=0$, proving the claim.
Assume then that there exists $g\in\pa (B,A)$ with $\beta :=gf\sim 1$,
and that $g=\{ g_n:B_n\rightarrow A_n\}$ is an inverse system of homomorphisms
commuting with the structure maps. By definition of homotopy, there exist
$r\ge 1$ and $\alpha_i\in \pa (A,A)$ with 
$1=\alpha^0\equiv\alpha^1\equiv\dots\equiv\alpha^r=\beta$. Because 
$\alpha:=\alpha^1\equiv 1$, for every $n\in\Bbb N$ there exists $m_0\ge n$
such that for $m\ge m_0$, $\tau_{mn} (\alpha *1)$ factors as follows:
$$
\CD
QA_m                   @>\alpha *1>>            A_m\\
@VVV                                            @VV\tau_{mn} V\\
\frac{QA_m}{qA_m^m}    @>h>>                    A_n\\
\endCD
$$
Therefore, given $a_1,\dots , a_m\in A_m$, we have:
$$
\multline
0=\tau (\alpha *1) (qa_1\dots qa_m)\\
 =\tau ((\alpha a_1-a_1)\dots (\alpha a_m -a_m))\\
 \equiv (-1)^m\tau (a_1\dots a_m) \mod{<\tau\alpha a_1,\dots \alpha a_m>}\\
\endmultline
$$
Thus if $a_1,\dots ,a_m\in \Ker(\tau\alpha)$, we have $\tau(a_1\dots a_m)=0$.
We have proven the following statement:
$$
\multline
(\forall n\ge 1) (\exists m_0\ge n) \text{ and for each } m\ge m_0\\
\text{ an } N=N_m\ge m  \text{ such that } (\Ker \tau_{mn}\alpha_n)^N=0\\
\endmultline\tag6
$$
We are going to show next that if $\alpha$ satisfies (6) and $\gamma\equiv\alpha$,
then $\gamma$ satisfies (6) too. It will follow that $\beta$ --and then
also $f$-- satisfies (6), whence $K^\infty=0$ as we had to prove. So assume
(6) holds for $\alpha$ and let $\gamma:A\rightarrow A$ with 
$\gamma\equiv\alpha$. Proceeding as above, we can find, for each $n$, an
$m_1\ge m_0\ge n$ such that if $m\ge m_1$, then
$$
0\equiv (-1)^m\tau\alpha (a_1\dots a_m) \mod{<\tau\gamma a_1,\dots \gamma a_m>}
$$
In particular $\tau_{mn} (\Ker\gamma_m)^m\subset \Ker\tau_{mn}\alpha$ whence
for $N$ as in (6) we have\linebreak
$(\Ker \tau\gamma_m)^{mN}=0$. \qed
\enddemo
\bigskip
\proclaim{Corollary 3.5} A pro-algebra $A$ is weak equivalent to zero iff
$A^\infty =0$.
\endproclaim
\demo{Proof} If $A\sim 0$ then $UA\defor 0$ is a deformation by 3.2-i) and 
3.4. Therefore $UA^\infty=0$, whence $A^\infty=0$. The converse is trivial.
\qed\enddemo
\bigskip
\remark{Remark 3.6} We can now see how far $\pa$ is from being a closed
model category. Indeed: by 3.5 above, if $0\rightarrow A$ factors as 
a weak equivalence followed by a fibration, then $A\sim 0$. On the other
hand, if $TV$ is a tensor algebra then clearly $TV^\infty\ne 0$, despite
the fact that the map $0\rightarrowtail TV$ has the LLP with respect to all
fibrations.\endremark
\bigskip
\subhead{4. Nil-homotopy v. Polinomial homotopy}\endsubhead 
\bigskip
4.0. We want to compare
our nil-homotopy relation with the more usual notion of homotopy defined
via polynomial homotopies, as used for example to define Karoubi-Villamayor
$K$-theory ([KV]). Given two homomorphisms
$f,g\in \Cal{PA}(A,B)$, we shall write $f\poleq g$ if there exists a 
homomorphism $h:A\rightarrow B[t]$, with values in the polynomial ring on the
commuting variable $t$, such that the following diagram commutes:
$$
\xymatrix{A\ar[d]_{(f,g)}\ar[r]^h& B[t]\ar[dl]^{
(\epsilon_0,\epsilon_1)}\\
B\times B &}
$$
Here $\epsilon_i$ stands for ``evaluation at $i$" ($i=0,1$). Note $\epsilon_1$
is defined even if $B$ is not unital, in which case $t\notin B[t]$; we set
$\epsilon_1(\sum_{i=0}^n a_i t^i)=\sum_{i=0}^n a_i$. Also note
that the map $(\epsilon_0,\epsilon_1)$ is a fibration; a natural linear
section is given by $(b_0,b_1)\rightarrow b_0+b_1t$.
We observe that $\poleq$ is a reflexive and symmetric relation, and that
if $f\poleq g$ then $fh\poleq gh$ (whenever the composition makes sense).
It follows that the equivalence relation $\polho$ generated by $\poleq$ 
is preserved by composition
on both sides. Thus $B[t]$ plays the r\^ole the free path space of a topological
space plays in ordinary topological homotopy. We showed in 1.2 above that nil-homotopy
can be described in terms of higher fold cylinders. Analogously, polynomial
homotopy (or simply pol-homotopy) can be defined in terms of higher free path spaces. Set $B^I=B[t]$,
and define $B^{I^n}$ inductively by the pull-back square:
$$
\CD
\quad B^{I^n}    @>\ \ \qquad>>                      \ \  \quad B^{I^{n-1}}\\
@VVV             @VV\epsilon_0^{n-1} V\\
B^I        @>>>  B\\
\endCD
$$
We write $\epsilon_0^n$ and $\epsilon_1^n$ for the composite maps
$B^{I^n}\rightarrow B^I\overset\epsilon_0\to\rightarrow B$ and
$B^{I^n}\rightarrow B^{I^{n-1}}\overset\epsilon_1^{n-1}\to\rightarrow B$.
Thus $(\epsilon_0^n,\epsilon_1^n):B^{I^n}\rightarrow B\times B$
is a fibration, and two maps $f_0, f_1:A\rightarrow B$ are $\polho$ iff 
there exist $n$ and $h:A\rightarrow B^{I^n}$ such that $h\epsilon_i^n=f_i$.
\smallskip
We write $\spa^{pol}$ for the (strong) polynomial homotopy category,
and call a map $f\in\pa(A,B)$ a polynomial equivalence if its class
$[f]^{pol}$ is an isomorphism in $\spa^{pol}$. A typical polynomial
equivalence is the projection $B=\bigoplus_{n=0}^\infty B_n\fib B_0$
of a graded algebra or pro-algebra onto the part of degree zero, which
is homotopy inverse to the inclusion $B_0\hookrightarrow B$. A homotopy between
the composite $B\fib B_0\hookrightarrow B$ and the identity map is given
by $h:B\rightarrow B[t]$, $h(b)=bt^{deg(b)}$. Projections 
of the form $\bigoplus_{n=0}^\infty B_n\overset{pol}\to\defor B_0$ 
shall be called {\it graded deformations}. For example power deformations
are graded, because power algebras are. We also
consider the category $\pupa$ having as objects those of $\pa$ and
as homomorphisms from $A$ to $B$ the homotopy classes $[UA,UB]^{pol}$.
The relation between nil-homotopy and polynomial homotopy is established
by the following:
\bigskip
\proclaim{Theorem 4.1}  
\roster
\item"{(i)}" The functor $U:\pa\rightarrow\paq$ carries pol-homotopic maps 
to pol-homotopic maps.
\medskip
\item"{(ii)}" If $f,g:A\rightarrow B$ are nil-homotopic and if $A$
is quasi-free, then they are also pol-homotopic.
\medskip
\item"{(iii)}" The functor $\pa\rightarrow [\pa]^{pol}$ is the localization
at the class of graded deformations, and the functor
$\gamma':\pa\rightarrow\pupa$ is 
the localization at the union of the classes of nil-deformations and graded-deformations.
There is a category equivalence $\hpaq^{pol}\approx\pupa$.
\endroster\endproclaim
\demo{Proof} (i) It suffices to show that if $f,g:A\rightarrow B\in\pa$
satisfy $f\poleq g$, then $Uf\poleq Ug$. Let $H:A\rightarrow B[t]$ be
a homotopy from $f$ to $g$. Then $H'=H\pi^A:UA\rightarrow B[t]$ is a
homotopy from $f\pi^A$ to $g\pi^A$ and $Uf, Ug$ are liftings of $f\pi^A, g\pi^A$
to $UB$. Hence by [CQ-2, Lemma 9.1], we have $Uf\poleq Ug$.
\smallskip
(ii)By Theorem 2.1, the map $CylA\rightarrow A$ is a power deformation
retraction, hence a graded deformation. It follows that $\gatc\poleq\gatu$,
and then $f\polho g$. 
\smallskip
(iii) The proof of the first assertion is analogous to the proof of the
first assertion of Theorem 3.2-ii). Next, we must show that $\gamma'$
inverts both nil and graded  deformations and is initial among
functors with such property. That $\gamma'$ inverts graded deformations
follows from (i), and that it inverts nil-deformations from (ii) and 3.2.
If $F:\pa\rightarrow\Cal C$ inverts both types of deformation, then
$\tilde F:\pupa\rightarrow \Cal C$, $A\mapsto FA$, 
$\upa(A,B)\owns[f]\mapsto F\pi^B Ff (F\pi^A)^{-1}$ stisfies $\tilde F\gamma'=F$
and is the only such functor.\qed\enddemo
\bigskip

\subhead{5. Derived Functors}\endsubhead
\smallskip
\definition{Notations 5.0} Recall from [Q] that if 
$F:\Cal M\rightarrow\Cal M'$ is a functor
between model categories, then the total (left) derived functor 
$\eLF:Ho\Cal{M}\rightarrow Ho\Cal M'$ is the (left) derived functor of
the composite $\Gamma ' F:\Cal C\rightarrow Ho\Cal M'$ with respect to
the localization $\Gamma:\Cal M @>>> Ho\Cal M$. Similarly, given a 
category $C$ together with a functor $\Gamma:C\rightarrow C'$, and
a functor $F:\pa\rightarrow C$, we may (and do) consider the total left 
and right
derived functors of $F$ with respect to $\Gamma$ and to 
$\gamma:\pa @>>>\upa$ and $\gamma':\pa @>>>\pupa$.
\enddefinition
\bigskip
\example{Motivation 5.1} The following proposition generalizes a common 
 procedure for deriving 
 functors. As a motivation, recall the way crystalline (or infinitesimal)
 cohomology is
 defined for commutative algebras of finite type over a field of 
 characteristic zero. Given an algebra $A$ one chooses a smooth $k$-algebra
 $R$ and an epimorphism $p:R\fib A$ and defines $H_{cris}^*A$ as the cohomology
 of the (commutative) de Rham pro-complex $\Omega_{R/I^\infty}$ where 
 $I=\Ker p$ (cf. [H], [I]). The essential step in proving that 
 $H_{cris}^*$
 is well-defined is the observation that if $A$ above is quasi-free, then
 $\hat R_I$ is an algebra of power series over $A$, and that (continuous) 
 $H_{dR}^*$ 
 satisfies the Poincar\'e Lemma: $H_{dR}(A)\cong H_{dR}^*A[[t]]$ ([H]). 
 Here $H_{dR}^*A[[t]]\overset{def}\to =H^*(\varprojlim\Omega^*(A[t]/<t^n>)$. 
 Actually Poincar\'e Lemma is derived from the stronger fact that 
 $\Omega^*(A)\overset\sim\to\rightarrow \Omega^*A[t]$ is a homotopy
 equivalence of pro-complexes ([H]).  A non-commutative analogue of this
 construction was given by Cuntz and Quillen in [CQ2]. They showed that
 the non-commutative de Rham pro-complex $XUA$ associated
 to an associative algebra $A$ has the homotopy type of the periodic cyclic
 complex $\theta\Omega(A)$. In the framework of this paper, we interpret these
 results
 as saying that crystalline and periodic (co)-homology are respectively the
 derived functors of commutative and of non-commutative de Rham cohomology
(see 5.4 below).
 The next proposition gives sufficient and necessary conditions so that
 when the construction above is applied to an arbitrary functor $F$, the
 result represents the left derived functor $\eLF$. We call this condition 
 the Poincar\'e
 condition because it resembles the Poincar\'e Lemma quoted above.
 In both the commutative and non-commutative cases, one uses the fact that, 
 in characteristic zero, de Rham cohomology is invariant under polynomial 
 equivalence. Thus the Poincar\'e condition is automatic (see 5.3). However
 there are Poincar\'e functors which are not pol-homotopy invariant.
 For instance the Grothendieck group $K_0$ is nil-invariant
 (and therefore represents its derived functor) despite the fact that in 
 general, $K_0(A[t])\ne K_0(A)$.\endexample
 \bigskip
 \proclaim{Theorem-Definition 5.2}{\rm (Poincar\'e Functors)}
 Let $F:\pa\rightarrow\Cal C$ and $\Gamma :\Cal C\rightarrow\Cal C'$
 be functors. The following are equivalent:
 \roster
 \item"{(i)}"  $FU$ represents the derived functor of $F$ with respect
 to $\Gamma$ and to $\gamma:\pa\rightarrow\upa$.
 \medskip
 \item"{(ii)}" $\Gamma FU$ is nil-invariant.
 \medskip
 \item"{(iii)}" Given any commutative diagram:
 $$
\xymatrix{R_0\ar[d]_{p_0}\ar[r]^{f} &R_1\ar[dl]^{p_1}\\
  A}
$$
where $p_i$ is a nil-deformation and $R_i$ is quasi-free ($i=0,1$),
 the map $\Gamma Ff$ is an isomorphism in $\Cal C'$.
 \item"{(iv)}" Given any pro-vector space $V$ and any quasi-free pro-algebra
 $R$, the map $\Gamma F(R\hookrightarrow P_R(V))$ is an isomorphism in $\Cal C'$.
 \item"{(v)}" Condition (iv) holds for $V=A$ and $R=UA$ ($A\in\pa$).
 \endroster
\medskip
 We call $F$ a {\rm Poincar\'e} functor if it satisfies the equivalent
 conditions above.
 \endproclaim
 \demo{Proof} We mimic the proof of the fact that a functor between model 
categories which preserves homotopy equivalences between cofibrant
objects admits a derived functor ([Qui 1.4.1]).
\smallskip
(i)$\iff$(ii) That (i)$\Rightarrow$(ii) is clear.
Assume now $\Gamma FU$ is nil-invariant, and let 
$\hat F:\upa\rightarrow\Cal C'$ be the induced functor. We have to prove that
$\hat F=\eLF$, i.e. that $\Gamma FU=\hat F\gamma$ is equipped with a natural
map $\alpha:\Gamma FU @>>> \Gamma F$ such that if 
$\hat G:\upa @>>> \Cal C'$ is another functor and 
$\beta: G:=\hat G\gamma @>>> F$ is a natural map then $\beta$ factors
uniquely through $\alpha$. Let $\alpha=\Gamma F(\pi^A):
\Gamma FUA @>>> \Gamma FA$ and set 
$\overline\beta=(\beta U)(G\pi^A)^{-1}:GA @>>> \Gamma FUA$. Then 
$\overline\beta$ satisfies $\beta=\alpha\overline\beta$ and is the
only such map.
\smallskip
(ii)$\Rightarrow$(iii) The map $f$ is a strong equivalence because each $p_i$
is a deformation and $R_i$ is quasi-free. Therefore $\Gamma FUf$ is an 
isomorphism. On the other hand we have $\pi^{R_1} Uf=f\pi^{R_0}$ where each 
$\pi^{R_i}$ is a deformation retraction; thus it is enough to show that
each $\Gamma F\pi^{R_i}$ is an isomorphism. But if 
$\iota_i:R_i\rightarrow UR_i$ is a right inverse for $\pi^{R_i}$, then 
$\iota_i\pi^{R_i}:UR_i\rightarrow UR_i$ is a nil-equivalence, whence the 
proof reduces
to showing that if $g:UB @>>> UB$ is a strong equivalence, then $\Gamma Fg$
is an isomorphism. We know by hypothesis that $\Gamma FUg$ is an isomorphism, 
and we have $\pi^{UB}FUg=g\pi^{UB}$. But $\Gamma F\pi^{UB}$ must be an isomorphism,
because $\Gamma FU\pi^{B}$ is, and both $\pi^{UB}$ and $U\pi^B$ have a right
inverse in common; namely the map induced by $T\rho:TB @>>> T^2B$.
\smallskip
(iii)$\Rightarrow$(iv) Let $r:P_R(V)\defor R$ be the projection map.
Then $r$ is a deformation and is a retraction of the canonical inclusion.
Thus (iv) is a particular case of (iii), with $R_0=A=R$ and $R_1=P_R(V)$.
\smallskip
(v) is logically weaker than (iv).
\smallskip
(v)$\Rightarrow$(ii) By virtue of Example 2.2, if (v) holds, then
$\Gamma F$ sends homotopy equivalences $UA \rightarrow UB$
to isomorphisms, whence $\Gamma FU$ sends weak nil-equivalences to 
isomorphisms.\qed\enddemo
\bigskip
\proclaim{Corollary 5.3} 
If $F$ preserves either nil-deformation
retractions or graded deformations, then it is Poincar\'e. In the latter
case $FU$ represents the left derived functor with respect to both 
$\gamma:\pa @>>>\upa$ and to $\gamma':\pa @>>>\pupa$.
\endproclaim

\demo{Proof} That $F$ is Poincar\'e means that its restriction to $\paq$
preserves nil homotopy (cf. 3.2). Such is the case if $F$ preserves either
nil-homotopy or, by 4.1-ii),
pol-homotopy of arbitrary pro-algebras. The same argument as in the proof
of the theorem shows that, in the latter case, $FU$ also represents the
derived functor with respect to $\gamma'$.\qed\enddemo
\bigskip
\proclaim{Corollary 5.4} Let 
$X:\pa\rightarrow \Cal{PS}:=${\bf ((}Pro-Supercomplexes{\bf ))} be the functor
which assigns to every pro-algebra $A$ the de Rham pro-super complex $XA$ of
2.4 above.
Let $\Gamma:\Cal{PS} @>>> Ho\Cal{PS}$ be the localization at the class of homotopy
equivalences and let $\gamma$ and $\gamma'$ be as above. If the ground
field $k$ has $char(k)=0$ then the functor $X$ is Poincar\'e (relative to 
$\Gamma$ and to $\gamma$), and its left derived functor with respect to 
both $\gamma$ and $\gamma'$ is 
represented by the periodic cyclic pro-complex $\theta\Omega$ of 2.4 above.
\endproclaim
\demo{Proof} In characteristic zero, the functor $X$ preserves polynomial 
homotopy (e.g. by [Kas], or by [CQ2\&3]), whence it is Poincar\'e  and 
$XU$ 
represents $LX$ (by 5.3). On the other hand, in any characteristic, 
$XUA$ is homotopy 
equivalent to $\theta\Omega UA$,
 because $UA$ is quasi-free (e.g. by [P]). In characteristic 
zero, by virtue of Goodwillie's theorem ([G1], [CQ2]), $\theta\Omega UA$
has the homotopy type of $\theta\Omega A$. Summing up, if $char(k)=0$ then
$FU\approx\theta\Omega$ represents $LX$.\qed\enddemo
\bigskip
\remark{Remark 5.5} In characteristic $p>0$, the lemma above fails to hold.
Indeed, if $X$ were Poincar\'e then --by 5.2--the homology of the periodic 
cyclic complex 

$CP(P_0(k))=Hom(Xk, XP_0(k))$ 
should be zero, which --as
a straightforward calculation shows-- it is not. See also Lemma 6.6 below. 
\endremark
\bigskip
\bigskip
\subhead{6. The derived functors of rational $K$-theory and Cyclic Homology}\endsubhead
\bigskip
The purpose of this section is to show that the functor which assigns
to every $\rat$-pro-algebra its rational $K$-theory space is (almost)
a Poincar\'e functor, and that its left derived functor is essentially
the fiber of the Chern character with values in negative cyclic homology.
See Theorem 6.2 below for a precise statement.
The proof of Theorem 6.2 has two main ingredients. The first ingredient 
is Goodwillie's isomorphism 
$$
K_*^\rat(A,I)\cong HN_*(A,I)\tag7
$$ 
between the relative rational $K$-group of
a nilpotent ideal and its analogue in negative cyclic homology [G2]. 
Actually Goodwillie's result is stated and proven for unital algebras;
we shall use an adaptation of this that holds for arbitrary pro-algebras,
which is obtained in 6.1 below. This adaptation says that 
the relative $K$-group of an infinitesimal deformation is isomorphic to
the corresponding negative cyclic homology group, and essentially reduces 
the question
of the Poincar\'eness of $K$ to that of $HN$. The second ingredient is the
calculation of relative $HN$ for a power deformation. This calculation is
carried out without any hypothesis on the characteristic of $k$ 
(Proposition 6.8).
\bigskip
\bigskip
\subhead{\smc 6.0. The derived functor of rational $K$-theory}\endsubhead
\smallskip
We use the following model for the rational
$K$-theory of a unital algebra or ring:
$$
K^{\rat}(A):=\rat_\infty BGlA
$$
Here $Gl$ is the general linear group, and $B$ denotes the simplicial
set associated to the category of $Gl$. Thus for us $K^\rat(A)$ is a simplicial
set; note that its homotopy groups are precisely Quillen's rational
$K$-groups. For general, non-necessarily unital algebras over the ground
field $k$ we set:
$$
K^\rat(A):=\hofiber(K^\rat(\tilde{A}) @>>>K^\rat(k))
$$
Thus in general $K^\rat(A)$ depends on $k$, and coincides with the usual
rational $K$-group if $A$ is unital or more generally if
it is excisive for $K^\rat$. Now we extend this definition to the case
of pro-algebras, by taking homotopy inverse limits, as follows. If 
$A=\{ A_\lambda:\lambda\in\Lambda\}$ we put:
$$
K^\rat(A):=\underset\Lambda\to{\ili}\ \ K^\rat(A_\lambda)
$$
Next we generalize Goodwillie's isomorphism to the pro-algebra case;
we assume throughout that $chark=0$.
Recall from [G2] that the isomorphism (7) is induced by a natural
Chern character $K_*^\rat(A) @>>> HN_*(A):=HN_*(A/k)$ which is defined for 
every unital algebra $A$. By [W] this character may be realized as a 
simplicial map $ch:K^\rat(A) @>>>SN(A)$, where $SN$ is constructed as 
follows. First truncate the total chain complex for negative cyclic homology
to obtain a complex $CN^t$ such that $H_n(CN^t)=HN_n(A)$ ($n\ge 1$)
and $H_n(CN^t)=0$ if $n\le 0$. Next define $SN$ as the result of applying the
Dold-Kan correspondence to $CN^t$. Hence $SN$ is a connected, fibrant simplicial
set with $\pi_nSN(A)=HN_nA$ ($n\ge 1$), and the isomorphism (7) says that
the map between fibers $K^\rat(A,I) @>>> SN(A,I)$ is a weak equivalence.
If now $A$ is any --non necessarily unital-- algebra, and $I\vartriangleleft A$
is a nilpotent ideal, then we have weak equivalences:
$$
K^\rat (A,I)\cong K^\rat(\tilde{A},I) @>\overset ch\to\sim>>SN(\tilde{A},I)
\cong SN(A,I)\tag8
$$
He have thus extended (7) to non-unital algebras. If now  
$A=\{ A_\lambda: \lambda\in\Lambda\}$ is a pro-algebra, we set 
$SN(A)=\underset\Lambda\to\ili\ \ SN(A_\lambda)$, and write $ch:K^\rat(A/k)@>>> SN(A)$
for the map induced by passage to $\ili$. As $\ili$ preserves fibers, fibrations
and weak equivalences of fibrant s. sets, (cf. [BK]) it follows that the weak
 equivalences (8) hold for arbitrary deformations and pro-algebras.  We have
 proven:
\bigskip
\proclaim{Lemma 6.1} With the notations and definitions of 6.0 above,
there is a natural map of fibrant simplicial sets $ch:K^\rat(A)@>>>SN(A)$
which is defined for all pro-algebras $A$, and coincides with Goodwillie's
character in the case of unital algebras. If $f:A\defor B$ is a deformation,
then the induced map $K^\rat(f)\approx SN(f)$ is a weak equivalence.\endproclaim

\demo{Proof} See the discussion above.\qed\enddemo
\bigskip
\bigskip
{\bf 6.1.1.} In particular the lemma above holds if $f$ is a power deformation of quasi-free 
pro-algebras, whence --by Theorem 5.2-iv)-- $K^\rat$ will be Poincar\'e iff $SN$ is. In the
next subsection we compute the homotopy groups of $SN(f)$  for power deformations
of quasi-free pro-algebras and show that these are all zero except for $\pi_1$, 
which is nonzero. Thus the simplicial set $SN'$ obtained from the complex
$CN$ by truncating in degree 2, so that $\pi_n(SN')=HN_n$ if $n\ge 2$ and
zero otherwise is a Poincar\'e functor; further, its derived functor is null-homotopic,
cf. 6.9 below. It follows that the $K$-theory space obtained by the same
process as above using the elementary group instead of the general linear
group is a Poincar\'e functor. Explicitly, the functor:
$$
KE^\rat(A):=\underset\Lambda\to\ili\ \ \hofiber(\rat_\infty(E\tilde{A}@>>>
\rat_\infty Ek)\tag9
$$
is Poincar\'e. 
\bigskip
\proclaim{Theorem 6.2}({\rm The derived functor of $K$-theory}) 
\smallskip
The functor $A\mapsto K^\rat(A)$ is not Poincar\'e.
However, the functor $A\mapsto KE^\rat (A)$ of (9) above is, and therefore
it has a left derived functor $LKE^\rat$. Set $LK^\rat_n(A):=\pi_nLKE^\rat$; 
then
there is an exact sequence:
$$
\split
\dots @>>> HN_{n+1}A@>>>LK_n^\rat(A) @>>>K_n^\rat(A)@>>>H_n(A)@>>>\\
\dots @>>> HN_3(A) @>>>LK_2^\rat(A)@>>>K_2^\rat(A)@>>>HN_2(A)\\
\endsplit
$$
\endproclaim
\demo{Proof} The first two assertions follow from the discussion above and
6.9 below.
To prove the third assertion consider the exact sequence of $K$-groups
associated with the universal deformation $\pi^A:UA\defor A$. Then 
$LK_n^\rat(A)=K_n^\rat(UA)$ ($n\ge 2$) (by 5.2) and $K_n(\pi^A)\cong HN_n(\pi^A)$
($n\ge 1$) (by 6.1). Because $UA$ is quasi-free, $HN_n(UA)=0$ for $n\ge 2$, and 
therefore 
$HN_n(\pi^A)\cong HN_{n+1}(A)$, for $n\ge 2$. This proves that the sequence 
is exact at $LK_2^\rat(A)$ and to the left. By the same argument, the natural
map $HN_2(A)\hookrightarrow HN_1(\pi^A)$ is injective, whence 
$K_2^\rat(A)@>>>K^\rat_1(\pi^A)$ factors through $ch_2$. It follows that the
sequence is exact also at $K_2^\rat(A)$, completing the proof.\qed\enddemo
\bigskip
\bigskip
\bigskip 
\subhead{\smc 6.3. The derived functor of negative cyclic homology}\endsubhead
\bigskip
The purpose of this subsection is to compute the homotopy type
of the relative space $SN(P_A(V) @>>> A)$ associated with
a power deformation retraction of a quasi-free pro-algebra $A$ over
a field. We do not make any assumptions with regards to $chark$.
The calculation uses two lemmas (6.4 and 6.6) which show the patologies
that appear in characteristic $p>0$. In particular, 6.6 gives a different
proof of the fact that the de Rham pro-complex $X$ is Poincar\'e iff
$chark=0$. In Lemma 6.4 we give a formula for the homotopy type of the $X$ pro-complex 
of a free product. 
Recall that if $A$ and $B$ are algebras, then there is an isomorphism of  
vector spaces:
$$
A*B=A\oplus B\oplus T(A\otimes B)\oplus T(B\otimes A)
\oplus T(A\otimes B)\otimes A\oplus T(B\otimes A)\otimes B
$$
In particular, the natural inclusion $T(A\otimes B)\hookrightarrow A*B$
is an algebra homomorphism. Putting this map together with the natural
inclusions $A\hookrightarrow  A*B$ and $B\hookrightarrow A*B$, we get map of 
super complexes:
$$
XA\oplus XB\oplus XT(A\otimes B)\overset\iota\to\hookrightarrow X(A*B)
$$
As all the maps in the above discussion are natural, all of this generalizes
immediately to the case of pro-algebras. The following lemma may be regarded
as a particular, easy case of [FT, 3.2.1]. We give an independent proof in 
this particular case.
\bigskip
\proclaim{Lemma 6.4}({\rm Compare [FT, 3.2.1]}) Let $A$, $B$ be pro-algebras.
 There exist a natural map of pro-mixed complexes:
$\pi:X(A*B)\rightarrow XA\oplus XB\oplus XT(A\otimes B)$ such that $\pi\iota=1$
and a natural homotopy $h:1\sim\iota\pi$.\endproclaim

\demo{Proof}  By naturality, we may assume $A$ and $B$ are algebras.
The map $A*B\rightarrow A\times B$, $a\mapsto (a,0)$, $b\mapsto (0,b)$
induces a retraction $XA*B\rightarrow XA\oplus XB$. Write 
$XA*B=XA\oplus XB\oplus Y$. Thus 
$Y_0=U\oplus V:=T(A\otimes B)\oplus T(B\otimes A)\oplus T(A\otimes B)\otimes A\oplus T(B\otimes A)\otimes B$.
where $U$ is the sum of the first two terms and $V$ is the sum of the last
two. Further, one checks that:
$$
Y_1\cong T(A\otimes B)dA\oplus T(B\otimes A)dB
\oplus\widetilde{T(A\otimes B)}\otimes A dB\oplus\widetilde{T(B\otimes A)}
\otimes BdA
\cong Y_0
$$
Consider the maps: $\alpha:U\rightarrow U$,
$x_0y_0\dots x_ny_n\mapsto y_nx_0\dots y_{n-1}x_n$, and 
$\mu:V\rightarrow U$,
$x_0y_0\dots x_ny_nx\mapsto xx_0y_0\dots x_ny_n$. Under the identifications
above, the map $\iota_1$ sends $x\in T(A\otimes B)\cong\Omega^1T(A\otimes B)_\natural$
onto $x+\alpha x\in U$. Define a mixed complex map 
$\pi:Y\rightarrow XT(A\otimes B)$, 
$\pi_0(u_0,u_1,v_0,v_1)=x_0+\mu y_0+\alpha x_1+\alpha\mu y_1$,
$\pi_1(u_0,u_1,v_0,v_1)=u_0$, $u_i\in U$, $v_i\in V$; $0$ denotes the 
alphabetical order, and $1$ denotes the inverse order. One checks that
$\pi\iota=1$. Further the map $h:Y_0\rightarrow Y_1$, 
$h(x_0,x_1,y_0,y_1)=(0,x_1+\mu y_1,y_0,y_1)$ verifies $\iota_1\pi_1=hb$
and  $\iota_0\pi_0=bh$.\qed\enddemo
\bigskip
\proclaim{Corollary 6.5}({\rm Compare [CQ-3, 7.3]}) If $chark=0$, then there
is a homotopy equivalence of supercomplexes:
$X(A*B)\approx XA\oplus XB$.
\endproclaim
\demo{Proof} Immediate from the well-known calculation of the cyclic homology
of a tensor algebra (e.g. [FT, 2.3.1]).\qed\enddemo
\bigskip
\bigskip
\proclaim{Lemma 6.6} Let $A$ be an algebra, $V$ a vector space and $P_A(V)$ the
power pro-algebra. Give $TV$ and $T(A\otimes TV)$ a gradation by setting
$deg(a)=0$ and $deg(v)=1$ ($a\in A$, $v\in V$). Then there exists a natural 
homotopy equivalence of pro-mixed complexes 
$$
XP_A(V)\approx XA\oplus\{X^{deg\le n}TV\oplus X^{deg\le n}
T(A\otimes TV): {n\ge1}\}
$$
\endproclaim
\demo{Proof} By definition the power pro-algebra $P_A(V)$ is graded, and
the gradation is given by the prescription of the lemma. This gradation is
reflected by the $X$-complex; we have a degree decomposition:
$$
C:=X(P_A(V))=\{\oplus_{i=0}^{2n} X^{deg=i}(P_{A}(V)_{n+1})\}
$$
We observe that
for $i\le n$ the direct summand subcomplexes corresponding to degree $i$ in the $X$ complex
of $P_{A}(V)_n=A*TV/<V>^n$ and of $A*TV$ are isomorphic. Further the 
pro-complex $D:=\{\oplus_{i\ge n+1}^{2n} C_{n+1}^{deg=i}\}$ is the zero 
pro-complex, as the structure maps $\tau_{n,2n}^D$ are all zero. Therefore
$C$ is isomorphic to the pro-complex $\{X^{deg\le n}(A*TV)\}$. Now the lemma
is immediate from 6.4.\qed\enddemo
\bigskip
\remark{Remark 6.7} As the homotopy equivalence in the lemma above is
natural, it extends automatically to pro-algebras. Since on the other
hand the Hochschild, cyclic and related homology groups of a tensor
algebra are well known, one could conceivably write down explicitly all
the relative pro-homology groups for the projection $P_A(V) @>>>A$ in
any characteristic.
In the next proposition we calculate the negative cyclic group for
the particular case when $A$ is an algebra and $V$ is a vector space.
Since in characteristic zero $HH_0(TV)\cong HH_1(TV)$, our calculation
can also be derived from [G2] in this particular case.\endremark
\bigskip
\proclaim{Proposition 6.8} Let $k$ be a field of characteristic $p\ge 0$, 
and let $A$ be a quasi-free
$k$-pro-algebra. If $V$ is a pro-vector space and $f:P_A(V)@>>> A$ is the
natural projection, then $SN(f)$ is an Eilenberg-Maclane space $E(\Upsilon(A,V),1)$,
where $\Upsilon(A,V)$ is an abelian group which depends functorially on $A$ and
$V$. Explicitly if $A$ is a quasi-free algebra and $V$ is a vector space,
then $\Upsilon(A,V)=\Pi_{n=0}^\infty(C_n\oplus\bigoplus_{r\ge 0}D_{n,r})$ is
the infinite product of the following co-invariant spaces:
$$
\multline
C_n=(T^nV)_{\Bbb Z/n} \text{\qquad and}\\ 
D_{n,r}=(\bigoplus_{i_1+\dots +i_r=n} A\otimes T^{i_1}V\otimes\dots\otimes A\otimes T^{i_r}V)_{\Bbb Z/r}\\
\endmultline
$$
Here $\Bbb Z/n$ and $\Bbb Z/r$ act by 
$v_1\otimes\dots\otimes v_n\mapsto v_n\otimes v_1\otimes\dots\otimes v_{n-1}$
and by 
$a_1\otimes x_1\otimes\dots\otimes a_r\otimes x_r\mapsto a_r\otimes x_r\otimes a_1\otimes 
x_1\otimes\dots\otimes a_{r-1}\otimes x_{r-1}$
\endproclaim
\demo{Proof} By the cofinality theorem for $\ili$ ([BK]), we may assume
$A$ and $V$ are indexed by $\Bbb N$. Thus for $n\ge 1$ we have an exact sequence:
$$
0 @>>> {\varprojlim}^1 HN_n(f_i) @>>> \pi_n(SN(f)) @>>>\varprojlim H_n(f_i) @>>> 0\tag10
$$
Since $P_A(V)$ is quasi-free, the inverse system $\{ HN_n(f_i): i\in \Bbb N\}$ is isomorphic
to the inverse system $\{ HN_n(X(f)): i\in\Bbb N\}$ (here $X$ is regarded as a 
mixed complex).
Thus both ends in the exact sequence above are zero for $n\ge 2$. Furthermore
$SN(f)$ is connected by definition; this concludes the proof of the first
assertion. Assume now $A$ is a quasi-free algebra and $V$ is a vector space.
It follows form  6.6 that we have an isomorphism of pro-vector spaces 
$$
\multline
\{ HN_1(f_n): n\in \Bbb N\}\cong\{\bigoplus_{i=0}^n T^i V_{\Bbb Z_i}\}\oplus\\
\bigoplus_{i=0}^n\bigoplus_{r\ge 0}\bigoplus_{j_1+\dots j_r=i} 
(A\otimes T^{i_1}V\otimes\dots\otimes A\otimes V^{i_r})_{\Bbb Z/r}: 
n\in \Bbb N\}\\
\endmultline\tag11
$$
As every map in the pro-vector space of the right hand of (11) is
a surjection, the $\varprojlim^1$ term in (10) is zero, and the second 
assertion of the proposition follows.\qed\enddemo
\bigskip
\bigskip
\proclaim{Corollary 6.9} The functor $A\mapsto SN(A)$ is not Poincar\'e,
regardless of the characteristic of $k$. The functor $A\mapsto SN'(A)$
of 6.1.1 above is Poincar\'e (in any characteristic) and its left derived
functor is null homotopic.\qed\endproclaim

\enddocument